\documentclass[a4paper,12pt,reqno]{amsart}
\usepackage{amssymb}
\usepackage{amsmath}

\def\curv{R^{\!\!\!^{^{\circ}}}\,}
\def\bu{{\hskip-1pt}_\star{\hskip-1pt}}
\def\i{\,\lrcorner\,}
\def\a{\alpha}
\def\b{\beta}

\def\vs{\vskip .6cm}

\def\la{\langle}
\def\ra{\rangle}
\def\.{\cdot}
\def\O{\Omega}
\def\n{\nabla}
\def\nb{\bar\nabla}
\def\l{\lambda}
\def\s{\sigma}
\def\t{\tilde}
\def\beq{\begin{equation}}
\def\eeq{\end{equation}}
\def\bea{\begin{eqnarray*}}
\def\eea{\end{eqnarray*}}
\def\beaa{\begin{eqnarray}}
\def\eeaa{\end{eqnarray}}
\def\ba{\begin{array}}
\def\ea{\end{array}}

\def\f{\varphi}

\def\o{\omega}

\def\L{\Lambda}

\def\bp{\begin{proof}}
\def\r{\end{proof}}

\def\D{\Delta}

\def \RM{\mathbb{R}}

\def \CM{\mathbb{C}}

\def\End{{\rm End}}
\def\XX{\chi}
\def\T{{\mathcal T}}

\def\d{{\delta}}

\def\Ric{\mathrm{Ric}}

\def\id{\mathrm{id}}
\def\be{\begin{equation}}

\def\ee{\end{equation}}

\def\pr{\rm{pr}}
\def\tr{\mathrm{tr}}
\def\Aut{\mathrm{Aut }}
\def\Cas{\mathrm{Cas}}

\def\Sym{\mathrm{Sym}}

\def\so{\mathfrak{so}}
\def\su{\mathfrak{su}}
\def\R{\mathbb{R}}

\def\SU{\mathrm{SU}}
\def\SO{\mathrm{SO}}

\def\U{\mathrm{U}}

\def\psp{\psi ^+}

\def\Sym{\mathrm{Sym}}
\def\scal{\mathrm{scal}}

\def\ai{A_{e_i}}

\newtheorem{ede}{Definition}[section]

\newtheorem{epr}[ede]{Proposition}

\newtheorem{ath}[ede]{Theorem}

\newtheorem{elem}[ede]{Lemma}

\newtheorem{ecor}[ede]{Corollary}

\title{Infinitesimal Einstein deformations of nearly K\"ahler metrics}

\author{Andrei Moroianu and Uwe Semmelmann}

\address{Andrei Moroianu \\ CMLS\\ {\'E}cole Polytechnique \\ UMR 7640 du CNRS
\\ 91128 Palaiseau \\ France}
\email{am@math.polytechnique.fr}

\address{Uwe Semmelmann\\ Mathematisches Institut, Universit{\"a}t zu
K{\"o}ln\\
Weyertal 86-90 D-50931 K{\"o}ln, Germany}
\email{uwe.semmelmann@math.uni-koeln.de}

\date{\today}

\begin{document}

\begin{abstract}
It is well known that every 6-dimensional strictly nearly K\"ahler manifold
$(M,g,J)$ is Einstein with positive scalar curvature
$\scal>0$. Moreover, one can show that  the space $E$ of
co-closed primitive $(1,1)$-forms on $M$ is stable under the Laplace
operator $\Delta$. Let $E(\l)$ denote the $\l$-eigenspace of the
restriction of $\Delta$ to $E$. If $M$ is compact, and has normalized
scalar curvature $\scal=30$,
we prove that the moduli space of infinitesimal Einstein deformations
of the nearly K\"ahler metric $g$
is naturally isomorphic to the direct sum $E(2)\oplus E(6)\oplus
E(12)$. From \cite{mns}, the last summand is itself isomorphic with
the moduli space of infinitesimal nearly K\"ahler
deformations.
\vs

\noindent
2000 {\it Mathematics Subject Classification}: Primary 58E30, 53C10, 53C15.

\medskip
\noindent{\it Keywords:} Einstein deformations, nearly
K\"ahler manifolds, Gray structures.
\end{abstract}

\maketitle

\section{Introduction}

Nearly K\"ahler manifolds, introduced by Alfred Gray in the 70s in the
framework of weak holonomy, are defined as almost Hermitian manifolds
$(M,g,J)$ which are not far from being K\"ahler in the sense that the
covariant derivative of $J$ with respect to the Levi-Civita connection
of $g$ is totally skew-symmetric.

The class of nearly K\"ahler manifolds is clearly stable under
Riemannian products. Using the generalization by Richard Cleyton and
Andrew Swann of
the Berger-Simons holonomy theorem to the case of connections with
torsion \cite{cs}, Paul-Andi Nagy showed in \cite{na} that every nearly
K\"ahler manifold is locally a Riemannian product of K\"ahler
manifolds, 3-symmetric spaces, twistor spaces over positive
quaternion-K\"ahler manifolds and 6-dimensional nearly K\"ahler
manifolds. This result shows, in particular, that genuine nearly
K\"ahler geometry
only occurs in dimension 6. It turns out that in this dimension,
strict ({\em i.e.} non-K\"ahler) nearly
K\"ahler manifolds have several other remarkable features: They carry
a real Killing spinor -- so they are in particular Einstein manifolds
with positive scalar curvature -- and
they have a $\SU_3$ structure whose intrinsic torsion is parallel
with respect to the minimal connection (cf. \cite{cs}). A strict
nearly K\"ahler structure on a compact
6-dimensional manifold with normalized scalar curvature $\scal=30$ is
called a Gray structure.

In \cite{mns} we have studied the moduli space $\mathcal G$
of infinitesimal deformations
of Gray structures on compact 6-dimensional manifolds, and
showed that this space is isomorphic to the space $E(12)$, where
$E(\l)$ denotes the intersection of the $\l$-eigenspace of
the Laplace operator and the space of co-closed primitive
$(1,1)$-forms.

In the present paper we consider the related problem of describing the
moduli space $\mathcal E$ of Einstein deformations of a Gray
structure. Since every Gray structure is in particular
Einstein, one has {\em a priori} $\mathcal E\supset \mathcal G$.  Our
main result (Theorem \ref{main}) gives a canonical
isomorphism between $\mathcal E$ and the direct sum $E(2)\oplus
E(6)\oplus
E(12)$.

The main idea is the following. The space of infinitesimal Einstein
deformation  on every compact manifold consists of trace-free
symmetric bilinear
tensors in a certain eigenspace of a second order elliptic operator called
the Lichnerowicz Laplacian $\Delta_L$. On a 6-dimensional nearly K\"ahler
manifold, one can decompose every infinitesimal Einstein deformation
$H$ (viewed as symmetric endomorphism) into its parts $h$ and $S$
commuting resp.
anti-commuting with $J$. Under the $\SU_3$ representation, the
space of symmetric endomorphisms commuting with $J$ is isomorphic to
the space of $(1,1)$-forms and that of symmetric endomorphisms
anti-commuting with $J$ is isomorphic to
the space of primitive $(2,1)+(1,2)$-forms and one may interpret the
eigenvalue equation for $\Delta_L $ in terms of the forms $\f$ and
$\s$ corresponding to $h$ and $S$.
The problem is that $\Delta_L$ does not commute with the isomorphisms
above, because $J$ is
not parallel with respect to the Levi-Civita connection. It is thus
natural to introduce a modified Lichnerowicz operator $\bar\Delta_L$,
corresponding
to the canonical Hermitian connection, better adapted to the nearly K\"ahler
setting. It turns out that the eigenvalue equation  for $\Delta_L $ translates,
{\em via} $\bar\Delta_L$, into a differential system for $\f$ and $\s$
involving the usual form Laplacian, which eventually yields the
claimed result.

\section{Preliminaries}

\subsection{Notation}
In this section we introduce our objects of study and derive several
lemmas which will be needed later. Here and henceforth, $(M^{2m},g,J)$
will denote an almost Hermitian manifold with tangent bundle $TM$,
cotangent bundle $T^*M$ and tensor bundle $\T M$. We denote as usual
by $\L^{(p,q)+(q,p)}M$ the projection of the complex bundle
$\L^{(p,q)}M$ onto the real bundle $\L^{p+q}M$. The
bundle of $g$-symmetric endomorphisms $\Sym M$ splits in a direct sum
$\Sym M=\Sym^+M\oplus\Sym^-M$, of symmetric endomorphisms commuting
resp. anti-commuting with $J$. The trace of every element in
$\Sym^-M$ is automatically 0, and $\Sym^+M$ decomposes further $\Sym^+
M=\Sym^+_{\,0}M\oplus\langle\id\rangle$ into
its trace-free part and multiples of the identity.

\subsection{Nearly K\"ahler manifolds} An almost
Hermitian manifold $(M^{2m},g,J)$ is called {\em nearly K\"ahler} if
\beq\label{1nk}(\n_XJ)(X)=0,\qquad\forall\ X\in TM,\eeq
where $\n$ denotes the Levi-Civita connection of $g$.
The canonical Hermitian connection $\nb$, defined by
\beq\label{der}\nb_XY:=\n_XY-\frac12J(\n_XJ)Y,\qquad\forall\
X\in TM,\ Y\in\XX(M),\eeq
is a $\U_m$ connection on $M$ ({\em i.e.} $\nb g =0$ and $\nb J=0$) with
torsion $\bar T_XY=-J(\n_XJ)Y$. A
fundamental observation, which -- although not explicitly stated -- goes
back to Gray, is the fact that $\bar\n\bar T=0$ on every nearly
K\"ahler manifold (see \cite{bm}).

We denote as usual the K\"ahler form of $M$
by $\o:=g(J.,.)$. The tensor $\psp:=\n\o$ is totally skew-symmetric by
(\ref{1nk}). Moreover, since $J^2=-\id$, it is easy to check that
$\psp(X,JY,JZ)=-\psp(X,Y,Z)$. In other words, $\psp$ is a form of type
$(3,0)+(0,3)$. Let us now assume that the dimension of $M$ is $2m=6$ and that
the nearly K\"ahler structure is strict, {\em i.e.} $(M,g,J)$ is not
K\"ahler. The form $\psp$ can be seen as the real
part of a $\nb$-parallel complex volume form on $M$, so $M$ carries an
$\SU_3$ structure whose minimal connection (cf. \cite{cs}) is exactly
$\nb$.

Let $A\in\L^1M\otimes\End M$ denote the tensor $A_X:=J(\n_XJ)=-\psp_{JX}$.
We will sometimes identify the endomorphism $A_X$ with the corresponding form
in $\L^{(2.0)+(0,2)}M$, {\em e.g.} in formula (\ref{a2}) below. By
definition, we have $\n_X=\nb_X+\frac12
A_X$ on $TM$. In fact this relation can be extended on the whole
tensor bundle, provided we use the right extension for $A_X$.

\subsection{The induced action} On a manifold $M$, every endomorphism
$A$ of $TM$ extends
as derivation to the tensor bundle $\T M$. In fact if we identify
$\End(T_xM) $ with $\mathfrak{gl}(n,\R)$ this is precisely the Lie algebra
action on the defining representation of $\T M$. We denote by $A\bu$ this
induced action. For example, we have
$$A\bu\tau=-\tau\circ A, \qquad A\bu f = A\circ f-f\circ A,\qquad
\hbox{for}\ \tau\in T^*M,\ f\in\End M.$$
If $(M,g,J)$ is almost Hermitian and $f\in\Sym^+M$, let $\t
f$ denote the associated $(1,1)$-form $\t f:=g(Jf.,.)$, so in particular
$\t\id=\o$, where $\o$ denotes the K\"ahler form of $M$. We compute,
for later use:
\beq\label{p1}(f\bu\o)(X,Y)=-\o(fX,Y)-\o(X,fY)=-g(JfX,Y)-g(JX,fY)=-2\t f(X,Y).
\eeq
A similar calculation shows that
\beq\label{p2}S\bu\o=0,\qquad\forall\ S\in \Sym^-M.\eeq
Notice that the
map $\End M\to \End(\T M)$, $A\mapsto A\bu$ is a Lie algebra
morphism, {\em i.e.}
$$[A,B]\bu=[A\bu,B\bu],\qquad\forall\ A,B\in\End M,$$
which can be expressed as
\beq\label{li}A\bu(B\bu T)=(A\bu B)\bu T+B\bu(A\bu T),\qquad\forall\
A,B\in\End M,\ T\in\T M.
\eeq
A convenient way of writing the induced action of $A\in\End M$ on a
$p$-form $u$ is
$$
A\bu u = -A^*(e_i) \wedge e_i \lrcorner u,
$$
where $A^*$ is the adjoint of $A$ and $\{e_i\}$ is a local orthonormal
basis of $TM$. Here, as well as in the remaining part of this paper, we
adopt the Einstein convention of summation on the repeated subscripts.

Notice that by (\ref{der}), the extensions of $\n$ and $\nb$ to the
tensor bundle $\T M$ are related by
\beq\label{nb}\nb_X=\n_X-\frac12(A_X)\bu.\eeq

\subsection{Algebraic results on nearly K\"ahler manifolds}
Assume that $(M^6,g,J)$ is a strict nearly K\"ahler manifold and that
the metric on $M$ is normalized such that $\scal=30$. From \cite[Theorem
5.2]{gray} it follows that for every unit vector
$X$, the endomorphism $\n_XJ$ (which vanishes on the 2-plane spanned by
$X$ and $JX$) defines
a complex structure on the orthogonal complement of that 2-plane. Then
the same holds for $A_X$ (because $A_X=-\n_{JX}J$).

The exterior bundle $\L^2M$ decomposes into irreducible
$\SU_3$ components as follows:
$$\L^2M\simeq \L^{(2,0)+(0,2)}M\oplus\L^{(1,1)}_0M\oplus \RM\omega.$$
The map $X\mapsto X\i\psp$ identifies the first summand with $TM$,
and $h\mapsto g(Jh.,.)$ defines an isomorphism between $\Sym_0^+M$
and the second summand.

Similarly, one can decompose $\L^3M$ into irreducible $\SU_3$ components
$$\L^3M\simeq \L^{(3,0)+(0,3)}M\oplus\L^{(2,1)+(1,2)}_0M\oplus
\L^1M\wedge\omega.$$
The first summand is a rank 2 trivial bundle spanned by $\psp$ and its
Hodge dual $*\psp$, and the isomorphism $S\mapsto S\bu\psp$ identifies
$\Sym^-M$ with the second summand.

If $\{e_i\}$ denotes a local orthonormal basis of $TM$,
it is straightforward to check the following formulas:
\beq\label{a1}A_{e_i}A_{e_i}(X)=-4X,\qquad\forall\ X\in TM.
\eeq
\beq\label{a2}A_{e_i}\wedge A_{e_i}=2\o ^2.
\eeq

\begin{elem} $(i)$ For every $X\in TM$, with corresponding 1-form
$X^\flat$ one has
\beq\label{a3}A_X\bu\psp=-2X^\flat\wedge\o.
\eeq
$(ii)$ If $S$ is a section of $\Sym^-M$, then the following formula
holds for every $X\in TM$:
\beq\label{a4}A_X\bu(S\bu\psp)=2(SX)^\flat\wedge\o.\eeq
\end{elem}
\bp  $(i)$ An easy computation shows:
\bea A_X\bu\psp&=&A_Xe_i\wedge e_i\i\psp=A_XJe_i\wedge Je_i\i\psp=
A_{Je_i}X\wedge\ai\\&=&\ai(JX)\wedge\ai=\frac12 JX\i(A_{e_i}\wedge
A_{e_i})\stackrel{(\ref{a2})}{=}-2X\wedge\o.
\eea
$(ii)$
The symmetric endomorphism $A_X\bu S=A_X\circ S-S\circ A_X$ commutes with
$J$ and is trace-free. Consequently, by Schur's Lemma (cf.~\cite{mns}
for a more detailed argument)
\beq\label{df}(A_X\bu S)\bu\psp=0.
\eeq
Notice that if $X^\flat$ is the 1-form corresponding to $X$
(which we usually identify with $X$), then $f\bu X^\flat=-(fX)^\flat$
for every symmetric endomorphism $f$.
We then compute:
\bea A_{X}\bu (S\bu\psp)&\stackrel{(\ref{li})}{=}&(A_{X}\bu S)\bu
\psp+S\bu(A_X\bu\psp)\stackrel{(\ref{a3}),(\ref{df})}{=}-2S\bu(X^\flat
\wedge\o)
\\&=&2(SX)^\flat\wedge\o-2X\wedge (S\bu\o)\stackrel{(\ref{p2})}{=}
2(SX)^\flat\wedge\o.
\eea
\r


\section{The curvature operator}

Let $(M^n, g)$ be a Riemannian manifold. The curvature
operator $\mathcal R : \Lambda^2M \rightarrow \Lambda^2M$
is defined by the equation
$g(\mathcal R (X\wedge Y), U \wedge V) = g(R_{X,Y}V,
U)$, for any  vector fields $X,Y,U,V$ on $M$, identified with the
corresponding 1-forms {\em via} the metric.
In a local orthonormal frame $\{e_i\}$ it can be written as
\beq\label{fo} \mathcal R(e_i \wedge e_j) = \frac12 R_{ijkl} e_l \wedge e_k =
-\frac12e_k \wedge R_{e_i, e_j}e_k.
\eeq
Using the identification of 2-vectors and (skew-symmetric) endomorphisms
given by $(X\wedge Y)(Z): = g(X, Z)Y - g(Y, Z)X$, formula (\ref{fo}) yields
$\mathcal R(X\wedge Y) (Z)=-R_{X, Y}Z$. Notice that a
manifold with curvature operator
$\mathcal R = c \id$ has Ricci curvature $c(n-1)$ and in
particular the curvature operator of the sphere is a positive
multiple of the identity.

Let $EM$ be the vector bundle associated to the bundle of orthonormal
frames {\em via} some representation $\pi :
SO(n) \rightarrow \Aut(E)$. Every orthogonal automorphism $f$ of $TM$
defines in a canonical way an automorphism of $EM$, denoted, by a
slight abuse of notation, $\pi(f)$. The differential of $\pi$ maps
skew-symmetric endomorphisms of $TM$ (or equivalently elements of
$\L^2M$) to endomorphisms of $EM$.
The Levi-Civita connection of $M$ induces a connection on $EM$ whose
curvature $R^E$ satisfies $R^E(X,Y)=\pi_*(R(X,Y))=-\pi_*(\mathcal
R(X\wedge Y)).$ Notice that $\pi_*(A)$ is exactly $A\bu$ in the
notation of Section 2.

We now define the curvature endomorphism
$q(R) \in \End (EM)$ as
\begin{equation}\label{qr-endo}
q(R) := -\frac12(e_i \wedge e_j)\bu\mathcal R(e_i
\wedge e_j)\bu.
\end{equation}
For example the curvature endomorphism $q(R)$ on
the form bundle $EM = \Lambda^pM$ satisfies
\beq\label{qr}
q(R) = (e_j \wedge e_i \lrcorner)\circ(  R_{e_i, e_j}e_k \wedge
e_k \lrcorner).
\eeq
In particular we have $q(R) = \Ric$ on 1-forms, and
$
q(R)=- \Ric\bu-2\mathcal R
$
on 2-forms. 

It is easy to check that the action of  $q(R)$ is compatible with
the identification of $\Lambda^2M $ with the space of skew-symmetric
endomorphisms (and, more generally, with all $\SO_n$ equivariant isomorphisms):
\begin{elem}
Let $\varphi \in \Lambda^2M$ be a $2$-form with associated
skew-symmetric endomorphism $A$, {\em i.e. }
$\varphi(Y, Z) = g(AY, Z)$ for any vector fields $Y, Z$. Then
$$
(q(R)\varphi)(Y,Z) = g((q(R)A)Y, Z).
$$
\end{elem}

We now return to the case of a 6-dimensional strict nearly
K\"ahler manifold $(M^6,g,J)$ with scalar curvature $\scal=30$.

Let $R$ be the curvature of the Levi-Civita connection and let $\bar
R$ be the curvature of the canonical Hermitian connection
$\bar \nabla$. The following relation between $R$ and $\bar R$
is implicitly contained in \cite{gray}.
\begin{elem}\label{gray}
For any tangent vectors $W, X, Y, Z$ one has
\bea
R_{W X Y Z} &=&\bar R_{W X Y Z} -\frac14g(Y,W)g(X,Z)
+\frac14g(X,Y)g(Z, W)
\\
&&+\frac34g(Y,JW)g(JX,Z) -\frac34g(Y,JX)g(JW,Z)
-\frac12g(X,JW)g(JY,Z).
\eea
\end{elem}
\bp
The stated formula follows using equation (3.1) and the polarization
of equation~(5.1) from \cite{gray}.
Note that there is a different sign convention for the curvature
tensors in \cite{gray}.
\r

The Ricci curvature of $\bar R$ satisfies $\overline \Ric = 4g$. This
follows from the formula above and the fact that
$(M^6,g,J)$ is Einstein with $\Ric = 5 g$.

Replacing $R$ by $\bar R$ in formula (\ref{qr}) yields a curvature
endomorphism $q(\bar R)$. It is easy to check that the curvature
operator with respect to $\bar\nabla$, denoted by $\bar {\mathcal R}$,
is a section of  $\Sym(\Lambda^{(1,1)}_0M)$ so we can express
$q(\bar R)=-\sum \a_i\bu \bar {\mathcal R}(\a_i)\bu$
for any orthonormal basis $\a_i$ of $\Lambda^{(1,1)}_0M$. Since
$\Lambda^{(1,1)}_0\RM^6\simeq \su_3$,
we see that $q(\bar R)$ preserves
all tensor bundles associated to $\SU_3$ representations. Moreover, a
straightforward computation using the fact that $\a\bu J=0$ and
$\a\bu\psp=0$ for every
$\a\in\su_3$ yields:
\beq\label{t1} g((q(\bar R)h)J.,.)=q(\bar R)\f, \qquad\hbox{and}\qquad
q(\bar R)(S\bu\psp)=(q(\bar R)S)\bu\psp,
\eeq
for every sections $h \in \Sym^+M$ and $ S \in \Sym^-M$, where
$\f$ denotes the $(1,1)$-form defined by $\f=g(hJ.,.)$.

The following lemma describes the difference $q(R) - q(\bar R)$.
It is an immediate consequence of the curvature formula in Lemma~\ref{gray}.
We will denote with  $\Cas  =  \frac12 (e_i\wedge e_j)\bu(e_i \wedge
e_j)\bu$ the Casimir operator of $\so(n)$ acting on the representation
$E$ and at the same time the corresponding endomorphism of $EM$.

\begin{elem}\label{difference}
The difference $ q(R)  - q(\bar R) \in
\End (EM) $ is given as
$$
q(R)  -  q(\bar R)  =  - \frac14 \Cas +
\frac38 (e_i\wedge e_j)\bu(Je_i \wedge Je_j)\bu
 - \frac18 (e_i\wedge Je_i)\bu(e_k \wedge Je_k)\bu
$$
\end{elem}

In the remaining part of this section we will apply Lemma~\ref{difference}
in order to compute $q(R)-q(\bar R)$ on certain spaces of endomorphisms
and forms.

On the bundle $\End M$ we define projections $\pr_\pm$ by $\pr_\pm(H) =
\frac12 (H\mp JHJ)$. Then $H = \pr_+(H) + \pr_-(H)$ is the
decomposition of the endomorphism $H$ in a part commuting
resp. anti-commuting with $J$. Using this notation we find
for the first sum in the above equation
$$
(e_i\wedge e_j)\bu(Je_i \wedge Je_j)\bu
  =
\left\{
\begin{array}{llll}
& -2 \id      & \mbox{on} & TM\\
& -8 \pr_-    & \mbox{on} & \Sym M \\
& -8 \pr_+    & \mbox{on} & \Lambda^2_0 M
\end{array}
\right.
$$
Indeed for any tangent vector $v \in TM$ we have
$$
(e_i\wedge e_j)\bu (J e_i \wedge Je_j)\bu v
=
(e_i\wedge e_j)\bu (g(Je_i, v)Je_j  -  g(Je_j, v)Je_i)
=
-2 (Jv \wedge e_j)\bu Je_j
=
 -2v.
$$
We now recall
that for every skew-symmetric endomorphism $A\in
\so(TM)\cong \Lambda^2M$ and for every $H\in \End M$ we have
$A\bu H = [A, H]$. Hence, we obtain
\bea (e_i\wedge e_j)\bu (J e_i \wedge Je_j)\bu H
&=&[(e_i \wedge e_j),[(J e_i \wedge Je_j),H]]\\
&=&(e_i\wedge e_j) (J e_i \wedge Je_j) H
+ H (e_i\wedge e_j) (J e_i \wedge Je_j)\\
&&-2(e_i\wedge e_j) H(J e_i \wedge Je_j)\\
&=&-4H-2(e_i\wedge e_j) H(J e_i \wedge Je_j).
\eea
It remains to compute the endomorphism $B:=(e_i\wedge e_j)  H (J e_i \wedge
Je_j)$. Applying it to a vector $Y$ and taking
the scalar product with a vector $Z$ we find
\bea
g(BY, Z)
&=&-g(H(J e_i \wedge Je_j)Y, (e_i\wedge e_j)Z)\\
&=&-g(H[ g(Je_i, Y)Je_j - g(Je_j,Y)Je_i ], [g(e_i, Z)e_j -g(e_j,Z)e_i])
\\&=&-2( g(JZ,Y)g(HJe_j,e_j) + g(HJZ,JY)).
\eea
Hence the sum $B=(e_i\wedge e_j)  H (J e_i \wedge Je_j)$ equals
$2JHJ$ if the endomorphism $H$ is symmetric and $-2JHJ$ if $H$ is
skew-symmetric and $\tr(HJ)=0$.

Next we compute the second sum of Lemma~\ref{difference} on tangent
vectors and endomorphisms. Since $2J = (e_i\wedge Je_i) $ we
immediately obtain
$$
(e_i\wedge Je_i)\bu(e_j \wedge Je_j)\bu
=
\left\{
\begin{array}{llll}
& -4 \id      & \mbox{on} & TM\\
& -16 \pr_-   & \mbox{on} & \End M
\end{array}
\right.
$$

We next determine the two sums of Lemma~\ref{difference}
on the space of 3-forms. Recall the type decomposition
$$
\Lambda^3M = \Lambda^{(3,0)+(0,3)}M \oplus \Lambda^{(2,1)+(1,2)}M,
$$
which coincides with the eigenspace decomposition of $(J\bu)^2$, with
eigenvalue $-9$ on the first and eigenvalue $-1$ on the second summand.
For any 3-form $\alpha$ we define a new 3-form $\hat \alpha$ by the formula
$\hat \alpha(X, Y, Z) = \alpha(JX, JY, Z) + \alpha(JX, Y, JZ) +
\alpha(X, JY, JZ)$. Then $(J\bu)^2 \alpha = -3 \alpha + 2\hat \alpha$
and the two components of $\Lambda^3M$ may also be characterized by
$$
\alpha\in \Lambda^{(3,0)+(0,3)}M     \qquad\mbox{if and only
  if}\qquad \hat \alpha =-3\alpha,
$$
and similarly
$$
\alpha \in\Lambda^{(2,1)+(1,2)}M     \qquad\mbox{if and only
  if}\qquad \hat \alpha =\alpha.
$$
Let $\pr_{3,0}$ and $\pr_{2,1}$ denote the projections onto the two
summands of $\Lambda^3M$. Then for $\alpha \in \Lambda^3M$
$$
\pr_{3,0}(\alpha) = \tfrac14(\alpha - \hat \alpha)
\qquad\mbox{and}\qquad
\pr_{2,1}(\alpha) = \tfrac14(3\alpha + \hat \alpha).
$$
For any 3-form $\alpha$ we have $(e_i \wedge Je_i)\bu (e_k \wedge
Je_k)\bu \alpha = 4 (J\bu)^2 \alpha = -12 \alpha + 8 \hat\alpha$,
which gives the second sum of Lemma~\ref{difference} and a
simple calculation yields  $(e_i \wedge e_j)\bu (Je_i \wedge
Je_j)\bu \alpha = -6\alpha -4\hat \alpha$ for the first sum.

Finally we may substitute the Casimir eigenvalues and our explicit
expressions into the formula of Lemma~\ref{difference}.
Recall that in the normalization with $\{e_i \wedge e_j\}$ as
ortho\-normal basis of $\Lambda^2 M\cong \so(TM)$, the Casimir
operator acts as $-p(n-p)\id$ on $\Lambda^pM$, and as
$-2n\id$ on $\Sym M$. Hence we obtain for $n=6$:

\begin{epr}\label{p33}
$$
q(R)  -  q(\bar R)   =
\left\{
\begin{array}{lrlll}
&\id  & & \hbox{\rm on} & TM\\
&  3 \pr_+ \!\!\!\!\! & + \, 2\, \pr_- & \hbox{on} & \Sym M \\
& - \pr_+  \!\!\!\!\! & + \,4 \,\pr_- & \hbox{on} & \Lambda^2_0 M\\
& -\pr_{2,1}  \!\!\!\!\! & +\,  9 \,\pr_{3,0}& \hbox{on} & \Lambda^3 M
\end{array}
\right.
$$
\end{epr}

\begin{ecor}\label{c34}
$$
q(R)  -  q(\bar R)   =
\left\{
\begin{array}{llll}
& \, 3\, \id & \mbox{on} & \Sym^+M \\
&  \,2 \,\id & \mbox{on} & \Sym^-M \\
& - \id & \mbox{on} & \Lambda^{(1,1)}_0 M\\
& - \id & \mbox{on} & \Lambda^{(2,1)+(1,2)}M
\end{array}
\right.
$$
\end{ecor}

\medskip


\section{Comparing rough Laplacians}

Let $(M,g,J)$ be a strict nearly K\"ahler manifold with Levi-Civita connection
$\n$ and canonical Hermitian connection $\nb$.
In this section we compare the actions of the rough Laplacians
$\n^*\n$ and $\nb^*\nb$ on several tensor bundles.

We will perform all calculations below at some fixed point $x\in M$
using a local orthonormal frame $\{e_i\}$ which is $\n$-parallel at
$x$. On any tensor bundle on $M$ we can write
$\n^*\n=-\n_{e_i}\n_{e_i}$
and because $\nb_{e_i}e_i=\n_{e_i}e_i-\frac12 J(\n_{e_i}J)e_i=0$, we
also have
$\nb^*\nb=-\nb_{e_i}\nb_{e_i}.$
We are interested in the operator $P:=\n^*\n-\nb^*\nb$. Using (\ref{nb})
and the fact that the tensor $A:=J\n J$ is $\nb$-parallel, we have
\bea P&=&-\n_{e_i}\n_{e_i}+\nb_{e_i}\nb_{e_i}=
-(\nb_{e_i}+\frac12 A_{e_i}\bu)(\nb_{e_i}+\frac12 A_{e_i}\bu)+
\nb_{e_i}\nb_{e_i}\\
&=&-\frac14 A_{e_i}\bu A_{e_i}\bu - A_{e_i}\bu\nb_{e_i}.
\eea

We now compute the action of the two operators occurring in the
previous formula on several tensor bundles which are
of interest in the deformation problem.

\begin{elem} \label{l41}
Let $\f$, $\s$, $h$ and $S$ be sections of $\L^{(1,1)}_0M$,
  $\L^{(2,1)+(1,2)}_0M$, $\Sym^+M$ and $\Sym^-M$ respectively. Then
\beq\label{e1} A_{e_i}\bu( A_{e_i}\bu \f)=-4\f.
\eeq
\beq\label{e2} A_{e_i}\bu (A_{e_i}\bu \s)= -4\s.
\eeq
\beq\label{e3} A_{e_i}\bu (A_{e_i}\bu h)= -12h.
\eeq
\beq\label{e4} A_{e_i}\bu (A_{e_i}\bu S)= -8S.
\eeq
\end{elem}

\bp Let $X\in TM$ be a tangent vector. Since $A_X$ and $\f$ are
2-forms of type
$(2,0)+(0,2)$ and $(1,1)$ respectively, we get
\beq\label{k}A_{X}(e_k,e_j)\f(e_k,e_j)=2\la A_X,\f\ra=0,
\qquad\forall\ X\in TM. \eeq
We then compute
\bea A_{e_i}\bu (A_{e_i}\bu \f)&=& A_{e_i}e_j\wedge
e_j\i(A_{e_i}e_k\wedge e_k\i\f)\\
&=&A_{e_i}e_j A_{e_i}(e_k,e_j)\wedge (e_k\i\f)-A_{e_i}e_j\wedge
A_{e_i}e_k \f(e_k,e_j)\\
&\stackrel{(\ref{k})}{=}&A_{e_i}^2(e_k)\wedge (e_k\i\f)-e_j\i(A_{e_i}\wedge
A_{e_i}e_k \f(e_k,e_j))\\
&\stackrel{(\ref{a1})}{=}&-4e_k\wedge(e_k\i\f)-\frac12e_j\i
e_k\i(A_i\wedge A_i \f(e_k,e_j))\\
&\stackrel{(\ref{a2})}{=}&-8\f-e_j\i
e_k\i(\o ^2\f(e_k,e_j))=-8\f-2e_j\i(Je_k\wedge \o \f(e_k,e_j))\\
&=&-8\f-2\o\f(e_k,Je_k)+2Je_k\wedge Je_j\f(e_k,e_j))\\
&=&-8\f+2e_k\wedge e_j\f(e_k,e_j))=-8\f+4\f=-4\f.
\eea
In order to prove (\ref{e2}), we express $\s$ as
$\s=S\bu\psp$ for some section $S$ of $\Sym^-M$.
Using (\ref{a4}) we obtain
\bea A_{e_i}\bu (A_{e_i}\bu \s)&=&2A_{e_i}\bu(Se_i\wedge\o)=
2Se_i\wedge(\ai\bu\o)=2Se_i\wedge
J\psp_{e_i}e_j\wedge Je_j\\
&=&-2Se_i\wedge\psp_{e_i}Je_j\wedge Je_j=4Se_i\wedge\psp_{e_i}=
-4S\bu\psp=-4\s.
\eea
Now, for every endomorphism $H$ of $TM$, we have
\beq\label{in}A_{e_i}\bu (A_{e_i}\bu H)=A_{e_i}^2H+HA_{e_i}^2-2\ai H
\ai\stackrel{(\ref{a1})}{=} -8H-2\ai H
\ai.\eeq
If $h\in\Sym^+M$, let $\f(.,.)=g(Jh.,.)$ be its associated
$(1,1)$-form. By (\ref{e1}) we have for every tangent vectors $X,Y$
\bea -4\f(X,Y)&=& A_{e_i}\bu( A_{e_i}\bu \f)(X,Y)\\
&=&\f(\ai ^2X,Y)+\f(X,\ai ^2Y)+2\f(\ai X,\ai Y)\\
&\stackrel{(\ref{a1})}{=}&-8\f(X,Y)+2g(hJ\ai X,\ai Y)\\
&=&-8\f(X,Y)+2g(J\ai h\ai X,Y),
\eea
whence $\ai h\ai=2h$. This, together with (\ref{in}), yields
(\ref{e3}). If $S\in \Sym^-M$, using the fact that $A_{JX}=A_X\circ
J=-J\circ A_X$ for every $X$, we can write
$$\ai S\ai=A_{Je_i}SA_{Je_i}=-\ai JSJ\ai=-\ai S\ai,$$
which together with (\ref{in}) yields
(\ref{e4}).
\r

\begin{elem} The following relations hold:
\beq\label{g0}e_i\i(\ai\bu \f)=0,\qquad\forall\ \f\in\L^{(1,1)}_0M.
\eeq
\beq\label{g1}(\ai\bu h)(e_i)=0,\qquad\forall\ h\in\Sym^+_{\,0}M.
\eeq
\beq\label{g2}e_i\wedge (\ai\bu\f)=0,\qquad\forall\ \f\in\L^{(1,1)}_0M.
\eeq
\beq\label{g3}e_i\i(A_i\bu(S\bu\psp))=0,\qquad\forall\ S\in\Sym^-M.
\eeq
\end{elem}
\bp
Simple application of the Schur Lemma, taking into account the
decomposition of the exterior bundles into irreducible components with
respect to the $\SU_3$ action.
\r

\begin{elem}\label{l43}
Let $\f$ and $S$ be sections of $\L^{(1,1)}_0M$ and
  $\Sym^-M$ respectively. If $h$ and $\s$ are defined as usual by
  $g(Jh.,.):=\f(.,.)$ and $\s:=S\bu\psp$, then
\beq\label{f4} A_{e_i}\bu \nb_{e_i} \f= -(J\d\f)\i\psp.
\eeq
\beq\label{f3} A_{e_i}\bu \nb_{e_i} \s= -2\d S\wedge\o.
\eeq
\beq\label{f1} (A_{e_i}\bu \nb_{e_i} h)\bu\psp=-2\d h\wedge\o-4d\f.
\eeq
\beq\label{f2} A_{e_i}\bu \nb_{e_i} S= (\d
S\i\psp+\d(S\bu\psp))\circ J.
\eeq
Here $\d$ denotes the co-differential on exterior forms and the
divergence operator whenever applied to symmetric endomorphisms.
\end{elem}

\bp Since $\ai$ anti-commutes with $J$ and $\nb_{e_i} \f$ is of type
$(1,1)$, it follows that $A_{e_i}\bu \nb_{e_i} \f$ is a form of type
$(2,0)+(0,2)$, so there exists a vector field $\a$ such that
$A_{e_i}\bu \nb_{e_i} \f=\a\i\psp$. In order to find $\a$, we use the
relation $(\a\i\psp)\wedge\psp=\a\wedge\o ^2$ (see \cite{mns}) and compute:
\bea \a\wedge\o ^2&=&(\a\i\psp)\wedge\psp=(A_{e_i}\bu \nb_{e_i}
\f)\wedge\psp\\
&=&A_{e_i}\bu ((\nb_{e_i}
\f)\wedge\psp)-\nb_{e_i}\f\wedge (\ai\bu\psp)\stackrel{(\ref{a3})}{=}
2\nb_{e_i}\f\wedge e_i\wedge\o\\
&=&-\nb_{e_i}\f\wedge(Je_i\i\o ^2)=-Je_i\i(\nb_{e_i}\f\wedge\o ^2)+
(Je_i\i\nb_{e_i}\f)\wedge\o ^2\\
&=&J(e_i\i\nb_{e_i}\f)\wedge\o
^2\stackrel{(\ref{g0})}{=}J(e_i\i\n_{e_i}\f)\wedge\o ^2=
-J\d\f\wedge\o ^2,
\eea
so $\a=-J\d\f$, thus proving (\ref{f4}).
Using the fact that $\psp$ is $\nb$-parallel, we get:
\bea A_{e_i}\bu \nb_{e_i}\s&=&A_{e_i}\bu (\nb_{e_i}S\bu\psp)
\stackrel{(\ref{li})}{=}(A_{e_i}\bu \nb_{e_i}S)\bu\psp+\nb_{e_i}S\bu
(A_{e_i}\bu\psp)\\
&\stackrel{(\ref{a3}),(\ref{df})}{=}&-2\nb_{e_i}S\bu(e_i\wedge\o)
\stackrel{(\ref{p2})}{=} 2(\nb_{e_i}S)e_i\wedge\o=-2\d S\wedge\o.
\eea
This proves (\ref{f3}). We next use (\ref{g1}) to write
\beq\label{u1}\d h=-(\n_{e_i}h)e_i=-(\nb_{e_i}h)e_i\eeq
whence
\bea (A_{e_i}\bu \nb_{e_i}
h)\bu\psp&\stackrel{(\ref{li})}{=}&A_{e_i}\bu ((\nb_{e_i} h)\bu\psp)-
(\nb_{e_i} h)\bu(A_{e_i}\bu\psp)\\
&\stackrel{(\ref{a3})}{=}&2(\nb_{e_i} h)\bu(e_i\wedge\o)=2(\nb_{e_i}
h)e_i\wedge \o+2e_i\wedge((\nb_{e_i} h)\bu\o)\\
&\stackrel{(\ref{p1})}{=}&-2\d h\wedge\o-4e_i\wedge\nb_{e_i}\f=-2\d
h\wedge\o-4d\f.
\eea

In order to check (\ref{f2}) we first compute
\bea \d(S\bu\psp)&=&-e_i\i\n_{e_i}(S\bu\psp)\stackrel{(\ref{g3})}{=}
-e_i\i\nb_{e_i}(S\bu\psp)=e_i\i\nb_{e_i}(S(e_j)\wedge\psp_{e_j})\\&=&
e_i\i((\nb_{e_i}S)e_j\wedge\psp_{e_j})=-g(\d
S,e_j)\psp_{e_j}-(\nb_{e_i}S)e_j\wedge(e_i\i\psp_{e_j})\\
&=&-\d S\i\psp+(e_i\i\psp_{e_j})\wedge(\nb_{e_i}S)e_j.
\eea
Let $B$ denote the endomorphism of $TM$ corresponding to the 2-form
$\d(S\bu\psp)+\d S\i\psp$. By the calculation above we get
\bea
B(X)&=&((e_i\i\psp_{e_j})\wedge(\nb_{e_i}S)e_j)(X)=
\psp(e_j,e_i,X)(\nb_{e_i}S)e_j-(\nb_{e_i}S)(e_j,X)(e_i\i\psp_{e_j})\\
&=&(\nb_{e_i}S)(\psp_{e_i}X)+\psp_{e_i}((\nb_{e_i}S)X)=(\nb_{e_i}S)(\ai
JX)+\ai(J(\nb_{e_i}S)X)\\
&=&(\nb_{e_i}S)(\ai JX)-\ai((\nb_{e_i}S)JX)=-(\ai\bu(\nb_{e_i}S))(JX).
\eea
Replacing $X$ by $JX$ yields (\ref{f2}).
\r
From Lemma \ref{l41} and Lemma \ref{l43} we infer directly
\begin{ecor}\label{co} Let $\f$ and $S$ be sections of $\L^{(1,1)}_0M$ and
  $\Sym^-M$ respectively. If $h$ and $\s$ are defined by
  $g(Jh.,.):=\f(.,.)$ and $\s:=S\bu\psp$, then
\beq\label{h1} (\n^*\n-\nb^*\nb)\f=\f+(J\d\f)\i\psp.
\eeq
\beq\label{h2} (\n^*\n-\nb^*\nb)\s=\s+2\d S\wedge\o.
\eeq
\beq\label{h3} (\n^*\n-\nb^*\nb)h=3h+s,\qquad\hbox{where}\
s\in\Sym^-M,\ \hbox{and}\ s\bu\psp=2\d h\wedge\o+4 d\f.
\eeq
\beq\label{h4} (\n^*\n-\nb^*\nb)S=2S-(\d S\i\psp+\d\s)\circ J.
\eeq
\end{ecor}
Finally, we obtain the invariance of the space of primitive co-closed
$(1,1)$-forms under the Laplace operator:
\begin{epr}\label{p45} If $\f$ is a co-closed section of $\L^{(1,1)}_0M$, then
the same holds for $\D \f$.
\end{epr}
\bp The 2-form $\D\f$ is clearly co-closed since $\f$ is co-closed.
Using (\ref{h1}), Corollary \ref{c34} and the classical Weitzenb\"ock formula
on 2-forms yield
$$\D\f=(\n^*\n+q(R))\f=(\nb^*\nb+q(\bar R))\f.$$
The last term is a section of $\L^{(1,1)}_0M$ since both $\nb$ and
$q(\bar R)$ preserve this space.
\r


\section{The moduli space of Einstein deformations}

We now have all the ingredients for the main result of this paper:

\begin{ath} \label{main}
Let $(M^6,g,J)$ be a Gray manifold.
Then the moduli space of infinitesimal Einstein deformations
of $g$ is isomorphic to the direct sum of the spaces of
primitive co-closed $(1,1)$-eigenforms of the Laplace operator for
the eigenvalues $2$, $6$ and $12$.
\end{ath}
\bp
Let $g$ be an Einstein metric with $\Ric = E g$.
From \cite{ab}, Theorem 12.30, the space of infinitesimal Einstein
deformations of $g$ is isomorphic to the set of symmetric
trace-free endomorphisms $H$ of $TM$ such that $\d H=0$ and
such that $\Delta_L H = 2 E H$, where $\Delta_L = \nabla^*\nabla +
q(R)$ is the so-called Lichnerowicz Laplacian $\Delta_L$.
Remark that $q(R)=2\curv +2 E \id$ in the notation of \cite{ab}.

In our present situation the Einstein constant equals $E=5$, so the space of
infinitesimal Einstein deformations of $g$ is isomorphic to the set
of $H \in \Sym M$ with $\d H=0=\tr H$ such that
\beq\label{ein} (\n^*\n + q(R))H=10H.
\eeq
Let $h:=\pr_+H$ and $S:=\pr_-H$ denote the projections of $H$ onto
$\Sym^\pm M$. We define the primitive $(1,1)$-form
$\f(.,.):=g(Jh.,.)$ and the $3$-form $\s:=S\bu\psp$. The key idea is
to express (\ref{ein}) in terms of an exterior differential system for
$\f$ and $\s$. Using
Corollary \ref{c34} and Corollary \ref{co}, (\ref{ein}) becomes
$$(\nb^*\nb + q(\bar R))(h+S)=10(h+S)-(3h+s)-(2S-(\d
S\i\psp+\d\s)\circ J)-3h-2S,
$$
where $s$ is the section of $\Sym^-M$ defined in the second part of
(\ref{h3}). Since the operator $(\nb^*\nb + q(\bar R))$ preserves the
decomposition $\Sym M=\Sym^+M\oplus \Sym^-M$, the previous equation
is equivalent to the system
\beq\label{s1}
\begin{cases} (\nb^*\nb + q(\bar R))h=4h+\d\s\circ J,\\
(\nb^*\nb + q(\bar R))S=6S-s,\\
\d S=0.
\end{cases}
\eeq
Taking the composition with
$J$ and using (\ref{t1}), the first equation of (\ref{s1}) becomes
\beq\label{s2}(\nb^*\nb + q(\bar R))\f=4\f-\d\s.
\eeq
Similarly, taking the action on $\psp$ and using (\ref{t1}) and the
definition of $s$, the second equation of (\ref{s1}) becomes
\beq\label{s3}(\nb^*\nb + q(\bar R))\s=6\s-2\d h\wedge\o-4d\f.
\eeq
Notice that $\d h=\d H-\d S=0$, which can also be seen by examining
the algebraic types in equation (\ref{s3}). From (\ref{u1})
we get
$$0=\d h= -(\nb_{e_i}h)e_i=(\nb_{e_i}J\f)e_i=J(\nb_{e_i}\f)e_i=-J\d
\f,$$
so finally the system (\ref{s1}) is equivalent to
\beq\label{s4}
\begin{cases} (\nb^*\nb + q(\bar R))\f=4\f-\d\s,\\
(\nb^*\nb + q(\bar R))\s=6\s-4d\f,\\
\d \f=0.
\end{cases}
\eeq
Using Corollary \ref{c34} together with the equations (\ref{h1}) and
(\ref{h2}) (keeping in mind that $\d S=0$ and $\d\f=0$) we get
the two identities
$(\nb^*\nb + q(\bar R))\f=(\n^*\n + q( R))\f$ and  $(\nb^*\nb + q(\bar
R))\s=(\n^*\n + q( R))\s$. Hence the classical Weitzenb\"ock
formula for the Laplace operator on forms implies that (\ref{s4}) is
equivalent to
\beq\label{s5}\begin{cases} \Delta\f=4\f-\d\s,\\
\Delta\s=6\s-4d\f,\\
\d \f=0.
\end{cases}
\eeq
\begin{elem}
Let $E(\lambda)$ be the $\lambda$-eigenspace of $\Delta$ restricted
to the space of co-closed primitive $(1,1)$-forms.
Then the space of solutions of the system $(\ref{s5})$ is isomorphic
to  the direct sum $E(2)\oplus E(6)\oplus E(12)$. The isomorphism
can be written explicitly as
$$(\f,\s)\stackrel{\Psi}{\mapsto} (8\f+\d\s,*d\s,2\f-\d\s)\in
E(2)\oplus E(6)\oplus
E(12)$$ with inverse
$$(\a,\b,\gamma)\in E(2)\oplus E(6)\oplus E(12)\stackrel{\Phi}{\mapsto}
\bigg(\frac{\a+\gamma}{10},
\frac{3d\a-5*d\b-2d\gamma}{30}\bigg).$$
\end{elem}
{\it Proof.} The first thing to check
is the fact that $\Phi$ and $\Psi$ take values in the right
spaces.

Let $(\a,\b,\gamma)\in E(2)\oplus E(6)\oplus E(12)$ and $(\f,\s):=
\Phi(\a,\b,\gamma)$.
From Lemma 4.4
in \cite{mns}, the exterior derivative of every co-closed primitive
$(1,1)$-form is a primitive $(2,1)+(1,2)$-form. Thus $\f\in \O
^{(1,1)}_0M$ and $\s\in \O ^{(2,1)+(1,2)}_0M$. A simple calculation
shows that $(\f,\s)$
satisfy the system (\ref{s5}).

Conversely, let $(\f,\s)$ be a solution of (\ref{s5}) and
$(\a,\b,\gamma):=\Psi(\f,\s)$. Clearly $\a,\ \b$ and $\gamma$ are
co-closed. From Proposition \ref{p45} and the first equation of
(\ref{s5}) we see that $\d\s$ is a section of $\L^{(1,1)}_0M$, so the
same holds for $\a$ and $\gamma$. The fact that $*d\s$ is
a section of the same bundle follows from  Lemma 4.3 in \cite{mns}. A
direct check shows that $\D\a=2\a$, $\D\b=6\b$ and $\D\gamma=12\gamma$.

Finally, it is straightforward to check that $\Phi$ and $\Psi$  are
inverse to each other. This proves the lemma and the theorem.
\r

In order to apply this result, one should try to compute the spectrum of the
Laplacian on 2-forms on some explicit compact nearly K\"ahler
6-dimensional manifolds. Besides the sphere $S^6$ -- which has no
infinitesimal Einstein deformations (cf. \cite{ab}) -- the only known
examples are 3-symmetric spaces $\CM\rm{P}^3=\SO_5/\U_2$,
$F(1,2)=\SU_3/U_1\times \U_1$ and $\SU_2\times \SU_2=\SU_2\times
\SU_2\times\SU_2 
/\Delta$.

Computations of the Laplace spectrum using the Peter-Weyl
theorem show 
that $E(2)$ and $E(6)$ vanish on each of these spaces. Moreover,
$E(12)$ vanishes on $\CM\rm{P}^3$ and on $\SU_2\times \SU_2$,
and is 8-dimensional on $F(1,2)$ (cf. \cite{au}). As a consequence of
these facts, we deduce:
\begin{enumerate}
\item every infinitesimal Einstein deformation of the 6-dimensional
  3-symmetric spaces is an infinitesimal Gray deformation
  (cf. \cite{mns});
\item the nearly K\"ahler structure on $\CM\rm{P}^3$ and on
  $\SU_2\times \SU_2$ is rigid;
\item there is an 8-dimensional space of infinitesimal deformations of
  the nearly K\"ahler structure on $F(1,2)$.
\end{enumerate}


\labelsep .5cm

\end{document}